\documentclass[11pt,a4paper,reqno]{amsart}
\title[Operators preserving real-rootedness]{On operators on polynomials 
preserving real-rootedness and the Neggers-Stanley Conjecture}

\author{Petter Bränd\'en}
\address{Matematik,
  Chalmers tekniska h\"ogskola och G\"oteborgs universitet,\linebreak
  S-412~96  G\"oteborg, Sweden}
\email{branden@math.chalmers.se}
\date{\today}

\usepackage[english]{babel}
\usepackage{amsmath}
\usepackage{amsthm,amssymb,latexsym}
\usepackage{t1enc}
\usepackage{epic}
\usepackage{eepic}
\usepackage{xypic}
\usepackage{mathrsfs}
\xyoption{all}

\newtheorem{proposition}{Proposition}
\newtheorem{lemma}[proposition]{Lemma}
\newtheorem{corollary}[proposition]{Corollary}
\newtheorem{theorem}[proposition]{Theorem}

\theoremstyle{definition}
\newtheorem{definition}[proposition]{Definition}

\newcommand{\N}{\mathbb{N}}
\newcommand{\om}{\omega}

\newcommand{\Pos}{\mathbb{P}}

\newcommand{\Four}{\mathcal{L}_{\phi}}
\newcommand{\E}{\mathcal{E}}
\newcommand{\R}{\mathbb{R}}

\newcommand{\inter}{\mathscr{I}}

\newcommand{\Al}{\mathscr{A}}

\newcommand{\e}{\ar@{-}}
\newcommand{\de}{\ar@{.}}

\def\ch#1,#2,{\binom#1#2}
\def\emm#1,{{\em #1}}
\def\ba#1,{\overline{#1}}

\def\gen#1,{\langle #1 \rangle}
\def\geno#1,{\langle #1 \rangle_{\infty}}

\def\newop#1{\expandafter\def\csname #1\endcsname{\mathop{\rm
#1}\nolimits}}

\newop{hp}
\newop{da}
\newop{ea}
\newop{lnfs}
\newop{MAJ}
\newop{row}
\newop{span}
\newop{min}
\newop{max}
\newop{deg}
\newop{des}
\newop{Con}
\newop{Ann}
\newop{dim}
\newop{Int}
\newop{Im}
\newop{ker}
\newop{mult}
\newop{rank}
\newop{sign}
\newop{sgn}
\newop{TrCl}
\newop{inf}

\usepackage{graphicx}

\begin{document}
\nocite{*}\bibliographystyle{plain}
\maketitle
\thispagestyle{empty} 
\begin{abstract}
We refine a technique used in a paper by Schur on real-rooted polynomials. 
This amounts to 
an extension of a theorem of Wagner on Hadamard products of Toeplitz 
matrices. We also apply our results to polynomials for which the 
Neggers-Stanley 
Conjecture is known to hold. More precisely, we settle 
interlacing properties for $E$-polynomials of series-parallel posets and 
column-strict labelled Ferrers posets.
\end{abstract}

\section{Introduction} 
Several polynomials associated to combinatorial structures 
are known to have real zeros. In most cases one can say more about the 
location of 
the zeros, than just that they are on the real axis. The 
{\em matching polynomial} 
of a graph is not only real-rooted, but it is known that the matching 
polynomial of the graph obtained by deleting a vertex of 
$G$ interlaces that of $G$ \cite{heilmann-lieb}. The same is true for 
the {\em characteristic polynomial} of graph (see e.g., \cite{godsil}). If 
$A$ is 
a nonnegative matrix and $A'$ is the matrix obtained by either deleting a row 
or a column, then Nijenhuis \cite{nijenhuis} showed that the 
{\em rook polynomial} of $A'$ interlaces that of $A$. 

The Neggers-Stanley Conjecture asserts that certain polynomials associated 
to posets, see Section \ref{negger}, have real zeros; see \cite{brenti,reinerwelker,wagner1} for 
the state of the art. For classes of posets for which the conjecture is known 
to hold we will exhibit explicit interlacing relationships.     

The first part of this paper is concerned with operators on polynomials which 
preserve real-rootedness. 
The following classical theorem is due to Schur \cite{schur}: 
\begin{theorem}[Schur]\label{basethm}
Let $f = a_0 +a_1x + \cdots +a_nx^n$ and $g = b_0 + b_1x + \cdots + b_mx^m$ 
be polynomials in $\R[x]$. Suppose that $f$ and $g$ have only real zeros 
and that the zeros of $g$ are all of the same sign. Then the polynomial
$$
f \odot g := \sum_k k!a_kb_kx^k,
$$
has only real zeros. If $a_0b_0 \neq 0$ then all the zeros of $f \odot g$ 
are distinct.
\end{theorem}
In this paper we will refine the technique used in Schur's proof of the 
theorem to extend a theorem of Wagner \cite[Theorem 0.3]{wagner2}.     
The {\em diamond product} of two polynomials $f$ and $g$ is the polynomial 
$$
f \Diamond g = \sum_{n \geq 0}\frac{f^{(n)}(x)}{n!}\frac{g^{(n)}(x)}{n!}
x^n(x+1)^n. 
$$
Brenti \cite{brenti} conjectured an equivalent form of Theorem 
\ref{wagnerthm} and Wagner proved it in \cite{wagner2}. 
\begin{theorem}[Wagner]\label{wagnerthm}
If $f,g \in \R[x]$ have all their zeros in the interval $[-1,0]$ then 
so does $f \Diamond g$. 
\end{theorem}  
This theorem has important consequences in combinatorics \cite{wagner1}, 
and it also has implications to the theory of total positivity \cite{wagner2}. 

In the second part of the paper we settle 
interlacing properties for $E$-polynomials of series-parallel posets and 
column-strict labelled Ferrers posets. 

We will implicitly use the fact that the zeros of a polynomial are continuous 
functions of the coefficients of the polynomial. In particular, the 
limit of real-rooted polynomials will again be real-rooted. For a treatment 
of these matters we refer the reader to \cite{marden}.  
\section{Sturm sequences and linear operators preserving real-rootedness}
Let $f$ and $g$ be real polynomials. 
We say that $f$ and $g$ {\em alternate} if $f$ and $g$ are real-rooted  
and either of the following conditions hold:
\begin{itemize}
\item[(A)] $\deg(g)=\deg(f)=d$ and 
$$
\alpha_1 \leq \beta_1 \leq \alpha_2 \leq \cdots \leq \beta_{d-1} \leq \alpha_d \leq \beta_d,
$$
where $\alpha_1 \leq \cdots \leq \alpha_d$ and 
$\beta_1 \leq \cdots \leq \beta_d$ are the zeros of $f$ and $g$ respectively 

\item[(B)] $\deg(f)=\deg(g)+1=d$ and 
$$
\alpha_1 \leq \beta_1 \leq \alpha_2 \leq \cdots \leq \beta_{d-1} \leq \alpha_d
$$
where $\alpha_1 \leq \cdots \leq \alpha_d$ and 
$\beta_1 \leq \cdots \leq \beta_{d-1}$ are the zeros of $f$ and $g$ 
respectively.
\end{itemize}
If all the inequalities above are strict then $f$ and $g$ are said to 
{\em strictly alternate}. Moreover, if $f$ and $g$ are as in (B) then 
we say that $g$ {\em interlaces} $f$, denoted $g \preceq f$. In the strict 
case we write $g \prec f$. 
If the leading coefficient of $f$ is positive we say that $f$ is 
{\em standard}. 

For $z \in \R$ let $T_z : \R[x] \rightarrow \R[x]$ be the 
translation operator defined by $T_z(f(x)) = f(x+z)$.  
For any linear operator $\phi : \R[x] \rightarrow \R[x]$ we define a 
linear transform $\Four : \R[x] \rightarrow \R[x,z]$ by 
\begin{eqnarray}\label{defour}
\Four(f) &:=& \phi(T_z(f)) \nonumber \\
         &=& \sum_n \phi(f^{(n)})(x) \frac {z^n} {n!} \\
         &=& \sum_n \frac {\phi(x^n)}{n!} f^{(n)}(z). \nonumber
\end{eqnarray}  

\begin{definition}\label{al}
Let $\phi : \R[x] \rightarrow \R[x]$ be a linear operator and let 
$f \in \R[x]$. If $\phi(f^{(n)})=0$ for all 
$n \in \N$, we let $d_\phi(f)=-\infty$. Otherwise let $d_\phi(f)$ be the 
smallest integer $d$ such that $\phi(f^{(n)})=0$ for all $n > d$.

The set $\Al^+(\phi)$ is defined as follows: If 
$d_\phi(f) = -\infty$, or   
$d_\phi(f)=0$ and $\phi(f)$ is standard real- and simple-rooted, then 
$f \in \Al^+(\phi)$. Moreover, $f \in \Al^+(\phi)$ if  
$d=d_\phi(f) \geq 1$ and all of the following conditions are satisfied:  
\begin{itemize}
\item[(i)] $\phi(f^{(i)})$ is standard for all $i$ and 
$\deg(\phi(f^{(i-1)}))= \deg(\phi(f^{(i)}))+1$ 
for $1 \leq i \leq d$,
\item[(ii)] $\phi(f)$ and $\phi(f')$ have no common real zero,
\item[(iii)] $\phi(f^{(d)}) \prec  \phi(f^{(d-1)})$,
\item[(iv)] for all $\xi \in \R$ the polynomial 
$\Four(f)(\xi,z)$ is real-rooted.
\end{itemize} 
Let $\Al^-(\phi):=\{ -f: f \in \Al^+(\phi)\}$ and $\Al(\phi):=\Al^-(\phi) \cup \Al^+(\phi)$. 
\end{definition} 
The following theorem is the basis for our analysis:
\begin{theorem}\label{huvud}
Let $\phi : \R[x] \rightarrow \R[x]$ be a linear operator. If 
$f \in \Al(\phi)$ then $\phi(f)$ is real- and simple-rooted and  
if $d_\phi(f) \geq 1$ we have 
$$
\phi(f^{(d)}) \prec \phi(f^{(d-1)}) \prec \cdots \prec \phi(f') 
\prec \phi(f).
$$
\end{theorem}
Before we give a proof of Theorem \ref{huvud} we will need a couple of  
lemmas. Note that  
$\frac{\partial}{\partial z}\Four(f) = \Four(f')$ so by Rolle's Theorem 
we know that $\Four(f')$ is real-rooted (in $z$) if $\Four(f)$ is. By 
Theorem \ref{huvud} it follows that $\Al(\phi)$ is closed under 
differentiation. 
A (generalised) {\em Sturm sequence} is a sequence $f_0, f_1, \ldots, f_n$ 
of standard polynomials such that $\deg(f_i)=i$ for 
$0 \leq i \leq n$ and 
\begin{equation}\label{sgn}
f_{i-1}(\theta)f_{i+1}(\theta) < 0,
\end{equation}
whenever $f_i(\theta)=0$ and $1 \leq i \leq n-1$. If $f$ is a standard 
polynomial with real simple zeros, we know from Rolle's 
Theorem that the sequence $\{f^{(i)}\}_i$ is a Sturm sequence. The 
following lemma is folklore.    
\begin{lemma}\label{sturm}
Let $f_0, f_1, \ldots, f_n$ be a sequence of standard polynomials 
with $\deg(f_i)=i$ for $0 \leq i \leq n$. Then the following statements are 
equivalent:
\begin{itemize}
\item[(i)] $f_0, f_1, \ldots, f_n$ is a Sturm sequence, 
\item[(ii)] $f_0 \prec f_1 \prec \cdots \prec f_n$.
\end{itemize}
\end{lemma}
The next lemma is of interest for real-rooted polynomials  
encountered in combinatorics. 
\begin{lemma}\label{lc}
Let $a_{m}x^{m} +a_{m+1}x^{m+1} + \cdots +a_{n}x^{n} \in \R[x]$ be 
real-rooted with  
$a_{m}a_{n} \neq 0$. Then the sequence $a_i$ is strictly log-concave, 
i.e.,      
$$
a_i^2 > a_{i-1}a_{i+1},   \ \  \ \ (m+1 \leq i \leq n-1).
$$
\end{lemma}
\begin{proof}
See Lemma 3 on page 337 of \cite{levin}.
\end{proof}
\begin{proof}[Proof of Theorem \ref{huvud}]
Let $f \in \Al^+(\phi)$. 
Clearly we may assume that $d=d_\phi(f)>1$.   
We claim that for $1 \leq n \leq d-1$:
\begin{equation}\label{lemeq}
\phi(f^{(n)})(\theta)=0 \ \ \Longrightarrow \ \ 
\phi(f^{(n-1)})(\theta)\phi(f^{(n+1)})(\theta) < 0. 
\end{equation}
If $1 \leq n \leq d-1$ and $\phi(f^{(n)})(\theta)=0$, then by condition 
(ii) and (iii) of Definition \ref{al} we have that there are integers 
$0 \leq \ell <n < k \leq d$ with 
$\phi(f^{(\ell)})(\theta)\phi(f^{(k)})(\theta)\neq 0$. 
By Lemma \ref{lc} and the real-rootedness of $\Four(f)(\theta,z)$ 
this verifies (\ref{lemeq}).  

If $\phi(f^{(d)})$ is a constant then 
$\{\phi(f^{(n)})\}_n$ is a Sturm sequence. Otherwise  
let $g = \phi(f^{(d)})$. Then, since $g' \prec g \prec \phi(f^{(d-1)})$, we 
have that \eqref{sgn} is satisfied everywhere in the sequence 
$\{g^{(n)}\}_n\cup \{\phi(f^{(n)})\}_n$. This proves the theorem by 
Lemma \ref{sturm}. 
\end{proof}
In order to make use of Theorem \ref{huvud} 
we will need further results on real-rootedness and interlacings of 
polynomials. There is a characterisation of alternating 
polynomials due to Obreschkoff and Dedieu. Obreschkoff proved the 
case of strictly alternating polynomials, see \cite[Satz 5.2]{obreschkoff}, 
and Dedieu \cite{dedieu} generalised it in the case 
$\deg(f)=\deg(g)$. But his proof also covers this slightly more general 
theorem: 
\begin{theorem}\label{obreschkoff2}
Let $f$ and $g$ be real polynomials. Then $f$ and $g$ alternate (strictly 
alternate) if and only if all polynomials in the space 
$$
\{ \alpha f +\beta g : \alpha, \beta \in \R \}
$$
are  real-rooted (real- and simple-rooted).
\end{theorem}
A direct consequence of Theorem \ref{obreschkoff2} is the following 
theorem, which the author has not seen previously in the literature.     
\begin{theorem}\label{altercor}
If $\phi : \R[x] \rightarrow \R[x]$ is a linear operator preserving 
real-rootedness, then $\phi(f)$ and $\phi(g)$ alternate if 
$f$ and $g$ alternate. Moreover, if $\phi$ preserves 
real- and simple-rootedness then $\phi(f)$ and $\phi(g)$ strictly 
alternate if $f$ and $g$ strictly alternate. 
\end{theorem}

\begin{proof}
The theorem is an immediate consequence of Theorem \ref{obreschkoff2} since 
the concept of alternating zeros is translated into a linear condition.
\end{proof} 
\begin{lemma}\label{sums}
Let $0 \neq h, f,g \in \R[x]$ be standard and real-rooted. If 
$h \prec f$ and $h \prec g$, then 
$h \prec \alpha f + \beta g$ for all $\alpha,\beta \geq 0$ not both equal 
to zero.
\end{lemma} 
Note that Lemma \ref{sums} also holds (by continuity arguments) when all 
instances of $\prec$ are replaced by $\preceq$ in Lemma \ref{sums}.
\begin{proof}
If $\theta$ is a zero of $h$ then  clearly $\alpha f + \beta g$ has the 
same sign as $f$ and $g$ at $\theta$. Since $\{h^{(i)}\}_i \cup \{f\}$ is 
a Sturm sequence by Lemma \ref{sturm}, so is  
$\{h^{(i)}\}_i \cup \{\alpha f + \beta g\}$. By Lemma \ref{sturm} again 
the proof follows.
\end{proof}
We will need two classical theorems on real-rootedness. The 
first theorem is essentially due to Hermite and Poulain and the second is 
due to Laguerre.
\begin{theorem}[Hermite, Poulain]\label{hermitepoulain}
Let $f(x) = a_0 + a_1x + \cdots + a_nx^n$ and $g$ be real-rooted. Then the 
polynomial 
$$
f(\frac {d}{dx})g:=a_0g(x) +a_1g'(x) + \cdots +a_ng^{(n)}(x)
$$
is real-rooted. Moreover, if $x^N \nmid f$ and $\deg(g) \geq N-1$
then any multiple zero of 
$f(\frac {d}{dx})g$ is a multiple zero of~ $g$.
\end{theorem}
\begin{proof}
The case $N=1$ is the Hermite-Poulain theorem. A proof can be found in any 
of the references \cite{levin,obreschkoff,schur}. For the general result it 
will suffice to 
prove that if $\deg(g) \neq 0$ then any multiple zero of $g'$ is a 
multiple zero of $g$. Let 
$$
g = c_0 + c_1(x-\theta) + \cdots + c_M(x-\theta)^M,
$$ 
where $c_M \neq 0$, $M >0$ and $(x-\theta)^2|g'$. Then $c_1=c_2=0$ and  $M>2$. 
If $c_0 = 0$ we are done and if $c_0 \neq 0$ we have by Lemma \ref{lc} 
that $0 = c_1^2 > c_0c_2 = 0$, which is a contradiction. 
\end{proof}
\begin{theorem}[Laguerre]\label{laguerre}
If $a_0 + a_1x + a_2x^2 + \cdots + a_nx^n$ is real-rooted then so is 
$$
a_0 + a_1x + \frac{a_2}{2!}x^2 + \cdots + \frac{a_n}{n!}x^n.
$$ 
\end{theorem}
\begin{proof}
Claim (ii) can be derived from (i) when applied to $x^n$, (see 
\cite{brenti}), or from Theorem \ref{basethm} as in \cite{levin,schur}. 
\end{proof}
We are now in a position to extend Theorem \ref{wagnerthm}.
\begin{theorem}\label{diamondtheorem}
Let $h$ be $[-1,0]$-rooted and let $f$ be real-rooted.
\begin{itemize} 
\item[(a)] Then $f \Diamond h$ is real-rooted, and if $g \preceq f$ then  
$$
g \Diamond h  \preceq f \Diamond h.
$$ 
\item[(b)]
If $h$ is $(-1,0)$- and simple-rooted and $f$ is simple-rooted  
then $f \Diamond h$ is simple-rooted and 
$$
g \Diamond h \prec f \Diamond h,
$$
for all $g \prec f$.
\end{itemize} 
\end{theorem}
\begin{proof}
First we assume that $\deg(h)>0$ and that $h$ is standard,   
$(-1,0)$-rooted and has simple zeros. Let 
$\phi : \R[x] \rightarrow \R[x]$ be the linear operator defined by 
$\phi(f) = f \Diamond h$.

We will show that $f \in \Al^+(\phi)$ if $f$ is standard real- and 
simple-rooted. Clearly we may assume that $\deg(f) = d \geq 1$. 
Condition (i) of 
Definition~ \ref{al} follows immediately from the definition of the 
diamond product. Now,  
$f^{(d-1)}= ax + b$, where $a,b \in \R$ and $a > 0$  so   
\begin{eqnarray*}
\phi(f^{(d)}) &=& ah \ \  \mbox{ and }\\
\phi(f^{(d-1)}) &=& (ax+b)h + ax(x+1)h',
\end{eqnarray*}
and since $h \preceq (ax+b)h$ and $h \preceq  x(x+1)h'$ we have 
by the discussion following Lemma \ref{sums} that 
$h \preceq \phi(f^{(d-1)})$. If 
$\theta$ is a common zero of $h$ and $\phi(f^{(d-1)})$, then 
$
\theta(\theta+1)h'(\theta) =0,
$ 
which is impossible since  
$\theta \in (-1,0)$ and $h'(\theta) \neq 0$. Thus 
$\phi(f^{(d)}) \prec \phi(f^{(d-1)})$, which verifies condition (iii) of 
Definition \ref{al}. Given $\xi \in \R$ we have 
\begin{eqnarray*}
\Four(f)(\xi,z) &=& \sum_n \frac {h^{(n)}(\xi)}{n!n!}\xi^n(\xi+1)^n 
               \frac {d^nf(\xi+z)}{dz^n}\\
       &=& H_\xi(\frac d {dz})f(\xi+z),
\end{eqnarray*}
where 
$$
H_\xi(x) = \sum_n \frac {h^{(n)}(\xi)}{n!n!}\{ \xi(\xi+1)x\} ^n. 
$$
By Theorem \ref{laguerre}  $H_\xi$ is real-rooted, which by Theorem 
\ref{hermitepoulain} verifies condition~ (iv). 

Suppose that $\xi$ is a common zero of $\phi(f')$ and $\phi(f)$. 
From the definition of the diamond product it follows that 
$\xi \notin \{0,-1\}$, so $x^2 \nmid H_\xi(x)$.
Since $\xi$ is supposed to be a common zero of $\phi(f')$ and $\phi(f)$ we 
have, by \eqref{defour}, that $0$ is a multiple zero of $\Four(f)(\xi,z)$. 
It follows from Theorem \ref{hermitepoulain} that $0$ is a multiple zero of 
$f(z+\xi)$, that is,  
$\xi$ is a multiple zero of $f$, contrary to assumption that $f$ is 
simple-rooted. This 
verifies condition (ii), and we can conclude that $f \in \Al^+(\phi)$.
Part (b) of the theorem now follows from Theorem \ref{altercor}.

If $h$ is merely $[-1,0]$-rooted and $f$ is real-rooted then 
we can find polynomials $h_n$ and $f_n$ whose limits are $h$ and $f$ 
respectively, such  
that $h_n$ and $f_n$ are real- and simple-rooted and $h_n$ is 
$(-1,0)$-rooted. Now, $f_n \Diamond h_n$ is real-rooted by the above 
and, by continuity, so is $f \Diamond g$. The proof now follows 
from Theorem \ref{altercor}.     
\end{proof}
There are many products on polynomials for which a similar 
proof applies. With minor changes in the above proof,   
Theorem \ref{diamondtheorem} also holds for the product 
$$
(f,g) \rightarrow \sum_{n \geq 0}\frac{f^{(n)}(x)g^{(n)}(x)}{n!}
x^n(x+1)^n.
$$

\section{Interlacing zeros and the Neggers-Stanley Conjecture}\label{negger}
Let $P$ be any finite poset of cardinality $p$. An injective function 
$\om : P \rightarrow \N$ is called a {\em labelling} of $P$ and $(P,\om)$ 
is a called a {\em labelled poset}.  A $(P,\om)$-partition 
with largest part $\leq n$ is a map $\sigma : P \rightarrow [n]$ such that
\begin{itemize}
\item $\sigma$ is order reversing, that is, if $x \leq y$ then 
      $\sigma(x) \geq \sigma(y)$, 
\item if $x < y$ and $\om(x) > \om(y)$ then $\sigma(x) > \sigma(y)$.
\end{itemize}
The number of $(P,\om)$-partitions with largest part $\leq n$ is denoted 
$\Omega(P,\om,n)$ and is easily seen to be a polynomial in $n$. Indeed, 
if we let $e_k(P,\om)$ be the number of surjective $(P,\om)$-partitions  
$\sigma : P \rightarrow [k]$, then by a simple counting argument we have:
\begin{equation}\label{omegaeq}
\Omega(P,\om,x) = \sum_{k=1}^{|P|}e_k(P,\om)\binom x k.
\end{equation}
The polynomial $\Omega(P,\om,x)$ is called the {\em order polynomial} of 
$(P,\om)$. The $E$-{\em polynomial} of $(P,\om)$ is the polynomial 
$$
E(P,\om) = \sum_{k=1}^pe_k(P,\om)x^k,
$$
so $E(P,\om)$ is the image of $\Omega(P,\om,x)$ under the invertible 
linear operator $\E : \R[x] \rightarrow \R[x]$ which takes 
$\binom x k$ to $x^k$. 
The Neggers-Stanley Conjecture asserts that 
the polynomial $E(P, \om)$ is real-rooted for all choices of 
$P$ and $\om$. The conjecture has been verified for series-parallel 
posets \cite{wagner1}, column-strict labelled Ferrers posets \cite{brenti} 
and for all labelled posets having at most seven elements.      

There are two operations on labelled posets under which $E$-polynomials
behave well. The first operation is the {\em ordinal sum}: 

Let $(P,\om)$ and 
$(Q,\nu)$ be two labelled posets. 
The {\em ordinal sum}, $P \oplus Q$, of $P$ and $Q$ is the poset with 
the disjoint union of $P$ and $Q$ as underlying set and with partial order 
defined by $x \leq y$ if either $x \leq_P y$, $x \leq_Q y$, or 
$x \in P, y \in Q$. For $i=0,1$ let 
$\om \oplus_i \nu$ be any labellings of 
$P \oplus Q$ such that 
\begin{itemize}
\item $(\om \oplus_0 \nu)(x) < (\om \oplus_0 \nu)(y)$ if 
$\om(x) < \om(y)$, $\nu(x) < \nu(y)$ or $x \in P, y \in Q$.  
\item $(\om \oplus_1 \nu)(x) < (\om \oplus_1 \nu)(y)$ if 
$\om(x) < \om(y)$, $\nu(x) < \nu(y)$ or $x \in Q, y \in P$.
\end{itemize}
The following result follows easily by combinatorial reasoning:
\begin{proposition}\label{brentiplus}
Let $(P,\om)$ and $(Q,\nu)$ be as above. Then 
$$
E(P\oplus Q, \om \oplus_1 \nu) = E(P,\om)E(Q,\nu)
$$
and 
$$
xE(P\oplus Q, \om \oplus_0 \nu)=(x+1)E(P,\om)E(Q,\nu),
$$
if $P$ and $Q$ are nonempty.
\end{proposition}
\begin{proof}
See \cite{brenti,wagner1}.
\end{proof}
The {\em disjoint union}, 
$P \sqcup  Q$, of $P$ and $Q$ is the poset on the disjoint union with 
$x < y$ in  $P \sqcup  Q$ if and only if $x <_P y$ or $x <_Q y$. 
Let $\om \sqcup \nu$ be any labelling of $P \sqcup  Q$ such that 
$$
(\om \sqcup \nu)(x) < (\om \sqcup \nu)(y),
$$
if $\om(x) < \om(y)$ or $\nu(x) < \nu(y)$. It is immediate by 
construction that
$$
\Omega(P \sqcup  Q,\om \sqcup \nu)=\Omega(P,\om)\Omega(Q,\nu)
$$
Here is where the diamond product comes in. 
Wagner \cite{wagner1} showed that the diamond product satisfies  
\begin{equation}\label{distruct}
f \Diamond g = \E(\E^{-1}(f)\E^{-1}(g)),
\end{equation}
which implies: 
\begin{equation}\label{asna}
E(P \sqcup  Q,\om \sqcup \nu) = E(P, \om) \Diamond E(Q, \nu),
\end{equation}
for all pairs of labelled posets $(P, \om)$ and $(Q, \nu)$.

If $P$ is nonempty and $x \in P$ we let $P \setminus x$ be the 
poset on $P \setminus \{ x \}$ with the order inherited by $P$. If 
$(P, \om)$ is labelled then $P \setminus x$ is labelled with 
the restriction of $\om$ to $P \setminus x$.   
By a slight abuse of notation we will write $(P \setminus x, \om)$ for 
this labelled poset. 
A series-parallel labelled poset 
$(S, \mu)$ is either the empty poset, a one element poset or 
\begin{itemize}
\item[(a)] $(S,\mu) = (P \oplus Q, \om \oplus_0 \nu)$,  
\item[(b)] $(S,\mu) = (P \oplus Q, \om \oplus_1 \nu)$ or 
\item[(c)] $(S,\mu) = (P \sqcup Q, \om \sqcup \nu)$ 
\end{itemize}
where $(P, \om)$ and $(Q, \nu)$ are series-parallel. Note that if 
$(S, \mu)$ is series-parallel then so is $(S\setminus x, \mu)$ for 
all $x \in S$. Let $\inter$ denote the class of 
all finite labelled posets $(S, \mu)$ such that $E(S,\mu)$ is real-rooted and 
$$
E(S \setminus x,\mu) \preceq E(S,\mu),
$$
for all $x \in S$. Note that the empty poset and the singleton posets are 
members of $\inter$ which by the following theorem gives that 
series-parallel posets are in $\inter$.  
\begin{theorem}\label{closed}
The class $\inter$ is closed under ordinal sum and disjoint union.
\end{theorem}
\begin{proof}
Suppose that $(P, \om), (Q,\nu) \in \inter$. \\
(a): Let $(S, \mu) =  (P \oplus Q, \om \oplus_0 \nu)$. Now, if $y \in P$ 
we have 
$$
(S\setminus y, \mu) = (P \setminus y \oplus Q, \om \oplus_0 \nu).
$$
If $|P|=1$ then by Proposition \ref{brentiplus} we have 
$E(S\setminus y, \mu)=E(Q, \nu)$ and $E(S, \mu) = (x+1)E(Q, \nu)$ so 
$E(S\setminus y, \mu) \preceq E(S, \mu) $. If $|P|>1$ then  
\begin{eqnarray*}
xE(S \setminus y, \mu) &=&  (x+1) E(P \setminus y, \om)E(Q, \nu) \\
   &\preceq& (x+1)E(P, \om)E(Q, \nu) \\     
   &=& xE(S, \mu), 
\end{eqnarray*}
which gives $E(S \setminus y, \mu) \preceq E(S, \mu)$. A similar argument 
applies to the case $y \in Q$. \\ 
(b): The case 
$(S, \mu) =(P \oplus Q, \om \oplus_0 \nu)$ follows as in (a). \\
(c): $(S,\mu) = (P \sqcup Q, \om \sqcup \nu)$. If $y \in P$ we have 
by (\ref{asna}) and Theorem~ \ref{diamondtheorem}:
\begin{eqnarray*}
E(S\setminus y, \mu) &=& E(P\setminus y  \sqcup Q, \om \sqcup \nu) \\
&=& E(P \setminus y , \om) \Diamond E(Q, \nu) \\
&\preceq& E(P, \om) \Diamond E(Q, \nu) \\
&=& E(S, \mu).
\end{eqnarray*}
This proves the theorem.
\end{proof}
In \cite{simion} Simion proved a special case of 
the following corollary. Namely the case when $S$ is a disjoint union of 
chains and $\mu$ is order-preserving. 
\begin{corollary}
If $(S,\mu)$ is series-parallel and $x \in S$ then 
$$
E(S\setminus x, \mu) \preceq E(S,\mu).
$$
\end{corollary}
Next we will analyse interlacings of $E$-polynomials of Ferrers posets. 
For undefined terminology in what follows we refer the reader to 
\cite[Chapter 7]{stanley2}. 
Let $\lambda = (\lambda_1 \geq \lambda_2 \geq \cdots \geq \lambda_{\ell}>0)$ 
be a partition. The {\em Ferrers poset} $P_{\lambda}$ is the poset 
$$
P_{\lambda} = \{(i,j) \in \Pos \times \Pos : 1\leq i \leq \ell, 1 \leq j \leq 
\lambda_i\},
$$  
ordered by the standard product ordering. A labelling $\om$ of 
$P_{\lambda}$ is {\em column strict} if $\om(i,j) >\om(i+1,j)$ and 
$\om(i,j) <\om(i,j+1)$ for all $(i,j) \in P_\lambda$. If $\om$ is 
a column strict labelling then any $(P_\lambda,\om)$-partition must 
necessarily be strictly decreasing in the $x$-direction and weakly 
decreasing in the $y$-direction. 
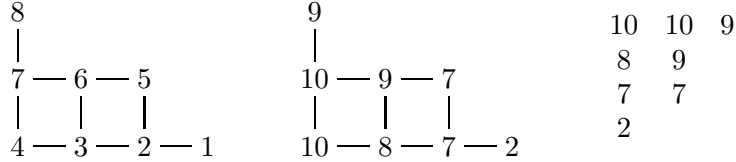
\begin{figure}\caption{From left to right: A column-strict labelling $\om$ of 
$P_\lambda$ with $\lambda=(3,2,2,1)$, a $(P_\lambda,\om)$-partition and 
the corresponding reverse SSYT. 
\label{fig2}}
$$
\vcenter{\xymatrix@R=12pt@C=12pt{
8\\
7\e[r]\e[u]&6\e[r]     &5          \\
4\e[r]\e[u]&3\e[r]\e[u]&2\e[r]\e[u]&1& \ \ 
    }}
\vcenter{\xymatrix@R=12pt@C=12pt{
9\\
10\e[r]\e[u]&9\e[r]     &7          \\
10\e[r]\e[u]&8\e[r]\e[u]&7\e[r]\e[u]&2& \ \ 
    }}
\begin{array}{ccc}
10&10&9\\
8&9\\      
7&7\\ 
2\\
\end{array}
$$
\end{figure}
It follows that the 
$(P_\lambda, \om)$-partitions are in a one-to-one correspondence with 
with the reverse SSYT's of shape $\lambda$ (see Figure \ref{fig2}). 
The number of reverse 
SSYT's of shape $\lambda$ with largest part $\leq n$ is by the combinatorial 
definition of the Schur function equal to $s_\lambda(1^n)$ which by 
the hook-content formula \cite[Corollary 7.21.4]{stanley2} gives 
us.  
\begin{equation}\label{hooklength}
\Omega(P_\lambda, \om, z)= \prod_{u \in P_\lambda}\frac{z+c_\lambda(u)}
{h_\lambda(u)},
\end{equation} 
where for $u=(x,y) \in P_\lambda$
$$
h_\lambda(u) := |\{ (x,j) \in \lambda : j \geq y\}|+ 
       |\{ (i,y) \in \lambda : i \geq x\}|-1
$$
and $c_\lambda(u) := y-x$  are the {\em hook length} respectively 
{\em content} at $u$. In \cite{brenti} 
Brenti showed that the $E$-polynomials of column strict labelled 
Ferrers posets are real-rooted. In the next theorem we refine this result. 
If $x<y$ in a poset $P$ and $x < z < y$ for no $z \in P$ we say that $y$ 
{\em covers} $x$. If we remove an element from $P_\lambda$ the resulting poset
will not necessarily be a Ferrers poset. But if we remove a maximal 
element $m$ from $P_\lambda$ we will have 
$P_\lambda \setminus m = P_\mu$ for a partition $\mu$ covered by 
$\lambda$ in the Young's lattice.
\begin{theorem}
Let $(P_\lambda, \om)$ be labelled column strict. Then $E(P_\lambda, \om)$ is 
real-rooted. Moreover, if $\lambda$ covers $\mu$ in the Young's lattice, then 
$$
E(P_\mu, \om) \preceq E(P_\lambda, \om).
$$
\end{theorem}
\begin{proof}The proof is by 
induction over $n$, where $\lambda \vdash n$. It is trivially true for 
$n=1$. If $\lambda \vdash n+1$ and $\lambda$ covers $\mu$ we have 
that $P_\lambda = P_\mu \cup \{m\}$ for some maximal element 
$m \in P_\lambda$. 
By definition $c_\mu(u)=c_\lambda(u)$ for all $u \in P_\mu$, so by 
(\ref{hooklength}) we have that for some $C > 0$:
$$
\Omega(P_\lambda, \om, x)=C(x+c_\lambda(m))\Omega(P_\mu, \om, x),
$$
and by \eqref{distruct}: 
$$
E(P_\lambda, \om) = C(x+c_\lambda(m))\Diamond E(P_\mu, \om).
$$
Wagner \cite{wagner1} showed that all real zeros of $E$-polynomials 
are necessarily in $[-1,0]$, so by induction we have that 
$E(P_\mu, \om)$ is $[-1,0]$-rooted. By  
Theorem~ \ref{diamondtheorem} this suffices to prove the theorem. 
\end{proof}

\bibliography{realroots}

\def\lfhook#1{\setbox0=\hbox{#1}{\ooalign{\hidewidth
  \lower1.5ex\hbox{'}\hidewidth\crcr\unhbox0}}}
\begin{thebibliography}{10}

\bibitem{brenti}
F.~Brenti.
\newblock Unimodal, log-concave and {P}\'olya frequency sequences in
  combinatorics.
\newblock {\em Mem. Amer. Math. Soc.}, 81(413):viii+106, 1989.

\bibitem{dedieu}
J.~Dedieu.
\newblock Obreschkoff's theorem revisited: what convex sets are contained in
  the set of hyperbolic polynomials?
\newblock {\em J. Pure Appl. Algebra}, 81(3):269--278, 1992.

\bibitem{godsil}
C.~D. Godsil.
\newblock {\em Algebraic combinatorics}.
\newblock Chapman and Hall Mathematics Series. Chapman \& Hall, New York, 1993.

\bibitem{heilmann-lieb}
O.~J. Heilmann and E.~H. Lieb.
\newblock Theory of monomer-dimer systems.
\newblock {\em Comm. Math. Phys.}, 25:190--232, 1972.

\bibitem{levin}
B.~Ja. Levin.
\newblock {\em Distribution of zeros of entire functions}.
\newblock American Mathematical Society, Providence, R.I., 1964.

\bibitem{marden}
M.~Marden.
\newblock {\em Geometry of polynomials}.
\newblock Second edition. Mathematical Surveys, No. 3. American Mathematical
  Society, Providence, R.I., 1966.

\bibitem{nijenhuis}
A.~Nijenhuis.
\newblock On permanents and the zeros of rook polynomials.
\newblock {\em J. Combinatorial Theory Ser. A}, 21(2):240--244, 1976.

\bibitem{obreschkoff}
N.~Obreschkoff.
\newblock {\em Verteilung und {B}erechnung der {N}ullstellen reeller
  {P}olynome}.
\newblock VEB Deutscher Verlag der Wissenschaften, Berlin, 1963.

\bibitem{reinerwelker}
V.~Reiner and V.~Welker.
\newblock On the {C}harney-{D}avis and the {N}eggers-{S}tanley conjectures.
\newblock {\em http://www.math.umn.edu/~reiner/Papers/papers.html}, 2002.

\bibitem{schur}
J.~Schur.
\newblock Zwei s\"atze \"{u}ber algebraische gleichungen mit lauter reellen
  wurzeln.
\newblock {\em J. Reine Angew. Math.}, 144(2):75--88, 1914.

\bibitem{simion}
R.~Simion.
\newblock A multi-indexed {S}turm sequence of polynomials and unimodality of
  certain combinatorial sequences.
\newblock {\em J. Combin. Theory Ser. A}, 36(1):15--22, 1984.

\bibitem{stanley2}
R.~P. Stanley.
\newblock {\em Enumerative combinatorics. {V}ol. 2}, volume~62 of {\em
  Cambridge Studies in Advanced Mathematics}.
\newblock Cambridge University Press, Cambridge, 1999.

\bibitem{wagner1}
D.~G. Wagner.
\newblock Enumeration of functions from posets to chains.
\newblock {\em European J. Combin.}, 13(4):313--324, 1992.

\bibitem{wagner2}
D.~G. Wagner.
\newblock Total positivity of {H}adamard products.
\newblock {\em J. Math. Anal. Appl.}, 163(2):459--483, 1992.

\end{thebibliography}

\end{document}